\documentclass[12pt,reqno]{amsart}
\usepackage{amssymb, amsmath,amsthm, multirow}
\usepackage{mathrsfs} 
\usepackage{bbm} 
\usepackage{mathabx} 
\usepackage{rotating} 
\usepackage{leftidx} 
\usepackage{hyperref}

\usepackage{stmaryrd} 
\usepackage{arydshln}
\makeatletter
\renewcommand*\env@matrix[1][*\c@MaxMatrixCols c]{%
  \hskip -\arraycolsep
  \let\@ifnextchar\new@ifnextchar
  \array{#1}}
\makeatother
\usepackage{pgf,tikz}
\include{master}
\usetikzlibrary{arrows, positioning, calc, chains}
\tikzset{
	ch/.style={circle,draw,on chain,inner sep=2pt},
	chj/.style={ch,join},
	every path/.style={shorten >=4pt,shorten <=4pt}
	}

\newtheorem{thm}{Theorem}[section]
\newtheorem{lem}[thm]{Lemma}
\newtheorem{prop}[thm]{Proposition}
\newtheorem{cor}[thm]{Corollary}

\newtheorem{lemma}[thm]{Lemma}

\theoremstyle{definition}

\newtheorem{exa}[thm]{Example}
\newtheorem{rmk}[thm]{Remark}

\newtheorem{alg}[thm]{Algorithm}

\theoremstyle{remark}

\numberwithin{equation}{section}
\newcommand{\uca}[1]{#1^{ \rotatebox{270}{$\curvearrowleft$}}}
\newcommand{\dca}[1]{#1^{ \rotatebox{270}{$\curvearrowright$}}}

\newcommand{\W}{W}
\newcommand{\tW}{\~{W}}
\newcommand{\D}{\mathscr{D}} 
\newcommand{\cH}{\mathbf{H}} 
\newcommand{\Sq}{\mathbf{S}} 

\newcommand{\ro}{\textup{ro}}
\newcommand{\co}{\textup{co}}

\newcommand{\lr}[1]{\left[ #1 \right]} 
\newcommand{\lrb}[2]{{#1 \brack #2}}

\newcommand{\dKq}{\dot{\mathbf{K}}}
\newcommand{\Uq}{\mathbf{U}_q}
\newcommand{\dUq}{\dot{\mathbf{U}}_q}

\newcommand{\g}[1]{\langle #1 \rangle}
\newcommand{\bp}[1]{\big{(} #1\big{)}}
\newcommand{\Bp}[1]{\Big{(} #1\Big{)}}

\newcommand{\ba}[1]{    \begin{array}   #1   \end{array}}

\newcommand{\bc}[1]{     \begin{cases}   #1   \end{cases}}
\newcommand{\bM}[1]{    \begin{bmatrix}    #1    \end{bmatrix}}

\renewcommand{\~}[1]{\widetilde{#1}}
\renewcommand{\^}[1]{\widehat{#1}}
\renewcommand{\=}[1]{\overline{#1}}

\newcommand{\ds}{\displaystyle}

\newcommand{\V}[1]{\widecheck{#1}}

\newcommand{\ZZ}{\mathbb{Z}}
\newcommand{\QQ}{\mathbb{Q}}

\newcommand{\NN}{\mathbb{N}}

\newcommand{\fa}{\textup{ for all }}
\newcommand{\fsm}{\textup{ for some }}

\newcommand{\rw}{\rightarrow}

\newcommand{\ld}{\lambda}
\newcommand{\Ld}{\Lambda}
\newcommand{\af}{\alpha}

\newcommand{\Tt}{\Theta}

\newcommand{\fgl}{\mathfrak{gl}}

\newcommand{\wo}{w_\circ}

\newcommand{\tif}{\textup{if }}

\newcommand{\Hom}{\textup{Hom}}

\newcommand{\otw}{\textup{otherwise}}

\newenvironment{enu}{\begin{enumerate}}{\end{enumerate}}
\usepackage{enumitem}
\newenvironment{enua}{\begin{enumerate}[label=$(\alph*)$]
}{\end{enumerate}}
\newenvironment{eq}{\begin{equation}}{\end{equation}}

\begin{document}
\title[Monomial bases of $\dot{U}_{\lowercase{q}}(\lowercase{\widehat{\mathfrak{gl}}_n})$]
{
    An elementary construction of monomial bases of modified quantum affine $\mathfrak{gl}_n$
}
\author{Chun-Ju Lai}
\address{
    Department of Mathematics, University of Virginia,
    Charlottesville, Virginia 22904;
    Current address: Max-Planck Institut f\"ur Mathematik in Bonn, Germany
}
\email{
    cl8ah@virginia.edu
}
\author{Li Luo}
\address{
	Department of mathematics, Shanghai Key Laboratory of Pure Mathematics and Mathematical Practice, East China Normal University, Shanghai 200241, China
}
\email{
    lluo@math.ecnu.edu.cn
}
\begin{abstract}
In 1990 Beilinson, Lusztig and MacPherson provided a geometric realization of modified quantum $\mathfrak{gl}_n$ and its canonical basis.
A key step of their work is a construction of a monomial basis.
Recently, Du and Fu provided an algebraic construction of the canonical basis for modified quantum affine $\mathfrak{gl}_n$,
which among other results used an earlier construction of monomial bases using Ringel-Hall algebra of the cyclic quiver.
In this paper, we give an elementary algebraic construction of a monomial basis for affine Schur algebras and modified quantum affine $\mathfrak{gl}_n$.

\end{abstract}
\maketitle
\date{}
\section{Introduction}
In \cite{BLM90}, Beilinson, Lusztig and MacPherson developed a geometric approach to the construction for the quantum group $\Uq(\fgl_n)$ together with its canonical basis.
Their work is based on a geometric realization of $q$-Schur algebras, which admits a standard basis, a monomial basis and a canonical basis.
Via a stabilization property for the family of $q$-Schur algebras, these bases naturally lift to the corresponding stabilization algebra, which is isomorphic to the modified quantum group (see also \cite{Lu93}).

The BLM construction has been partially generalized to quantum affine $\fgl_n$ by Ginzburg-Vasserot \cite{GV93} and  Lusztig \cite{Lu99}.
Recently, Du and Fu \cite{DF14} provided an algebraic BLM-type construction for quantum affine $\fgl_n$ via a Hecke-algebraic realization of $q$-Schur algebras.

In the original work, an essential step \cite[Proposition~3.8]{BLM90} towards proving the stabilization property is a direct construction of a  monomial basis for $q$-Schur algebra of finite type $A$ with some favorable properties.
By construction, each basis element is obtained by multiplying certain elements corresponding to divided powers of Chevalley generators.
However, since the Chevalley generators alone do not generate the entire quantum affine $\fgl_n$,
the construction above does not generalize easily to the $q$-Schur algebra $\Sq$ of affine type $A$.
Deng, Du and Xiao constructed a family of monomial bases for the Ringel-Hall algebra of the cyclic quiver using a ``strong monomial basis property'' (cf. \cite{DDX07}) for its composition subalgebra;
the construction therein is rather involved and it is not clear how to understand their bases along the line of \cite{BLM90}.
Using the surjection from the double Hall algebra of the cyclic quiver to $\Sq$ (proved in \cite{DDF12}), any such basis can then be adapted to a monomial basis for $\Sq$,
and this is the approach adapted in \cite{DF14}.

The goal of this paper is to provide a new and elementary construction (Theorem~\ref{thm:min gen}) of a monomial basis for $\Sq$
with similar favorable properties as in \cite[Proposition~3.8]{BLM90}.
Here each basis element is obtained by multiplying certain bar-invariant elements corresponding to bidiagonal matrices.
It turns out that these basis elements can be identified \cite{DDPW08} with the semisimple representations of Hall algebra of the cyclic quiver,
though this fact is not needed in our work.
Also, a variant of our construction gives a monomial basis for the stabilization algebra $\dKq$ that is isomorphic to the modified quantum affine $\fgl_n$.
As usual, the construction of canonical bases for $q$-Schur algebras and stabilization algebras follow readily from the existence of the corresponding monomial bases. Du and Fu \cite{DF14} provided such a construction of canonical bases for modified quantum affine $\mathfrak{gl}_n$.
The validity of our approach relies on the remarkable multiplication formulas \cite[Proposition~3.6]{DF15}, which are indeed the foundation of the BLM-type construction for quantum affine $\fgl_n$ and its canonical basis.

The BLM construction has also been generalized to finite types $B/C$ in \cite{BKLW14} and to finite type $D$ in \cite{FL14},
and it is natural to ask for a suitable affinization of these BLM-type constructions.
Our construction (without use of Hall algebras) was motivated by and in turn can be adapted to a construction of monomial bases of $q$-Schur algebras and stabilization algebras for affine type C (cf. a joint work \cite{FLLLW16} with Fan, Li and Wang), where no Hall-algebraic approach is yet available.

In Section 2, which is mostly expository, we review some necessary background and notations. 
 In Section 3 we construct a monomial basis for the affine $q$-Schur algebra (cf. Algorithm~\ref{alg:mono}) and show that it has some favorable properties (cf. Theorem~\ref{thm:min gen} and Corollary~\ref{cor:mA}).

Throughout this article, let $\NN=\{ 0,1,2,\ldots \}$ be the set of natural numbers.
Denote by $[a..b], [a..b), (a..b]$ and $(a..b)$ the integer intervals for $a,b\in\ZZ$. 
Let $v$ be an indeterminate over $\QQ$, and let $\lr{a}=\frac{v^{2a}-1}{v^2-1}$ for $a\in\ZZ$.
Finally, let $n, d$ be fixed positive integers.
\subsection*{Acknowledgments}
We thank Zhaobing Fan, Yiqiang Li and Weiqiang Wang for the ongoing collaboration,
which motivated this work. We thank Huanchen Bao and Qiang Fu for very helpful comments and suggestions.

Part of this work was done while the second author visited University of Virginia and when both of us visited Institute of Mathematics, Academia Sinica, Taipei.
We thank Shun-Jen Cheng for offering a wonderful working environment.
The second author is supported by Science and Technology Commission of Shanghai Municipality (Grant No. 13dz2260400) and National Natural Science Foundation of China (Grant No. 11671108).
\section{Affine Hecke algebras and affine $q$-Schur algebras}

\subsection{(Extended) affine Weyl groups}
Let $W$ be the Weyl group of type $\~{A}_{d-1}$
generated by $S=\{s_1, s_2, \ldots, s_d=s_0\}$.
The extended Weyl group  $\tW$ is generated by $W$ and $\pi$ satisfying $\pi s_i \pi^{-1}=s_{i-1}$ for $i=1,\ldots, d$.
It is well-known that $\tW$ can be identified as a permutation subgroup of $\ZZ$ satisfying
$g(i+d) = g(i)+d$ for all $i\in\ZZ, g\in\tW$.
In this identification each $s_i$ is mapped to the permutation $\prod_{k\in \ZZ} (kd+i , kd+i+1)$ and $\pi$ is mapped to the permutation $t\mapsto t+1$ for $t\in\ZZ$.
Denote the length function on $\W$ relative to $S$ by $\ell$.
Notice that each $g \in \tW$ can be uniquely expressed as $g=\pi^z w$
for some $z\in\ZZ$ and $w\in \W$, so the notion of length on $\W$
can be extended to $\tW$ by requiring $\ell(\pi) = 0$, or
equivalently, $\ell(g) = \ell(w)$.

The following lemma can be found in \cite[(3.2.1.1)]{DDF12}.
\begin{lemma}\label{lem:l(g)}
The length of $g \in \tW$ is given by
\[
\ell(g)
=
\sharp
\left\{
(i,j)\in[1..d]\times \mathbb{Z}
~|~
i<j, g(i)>g(j)
\right\}.
\]
\end{lemma}

Denote the set of (weak) compositions of $d$ into $n$ parts by
\[
\Ld = \Ld(n,d) = \{ \ld = (\ld_1, \ldots, \ld_{n})\in\NN^{n} ~|~ \textstyle\sum_{i=1}^n \ld_i = d \}.
\]
For each $\ld\in\Ld$, denote by $\W_\ld$ the parabolic subgroup of $\W$ with respect to $\ld$ generated by $S\backslash\{ s_{\ld_1}, s_{\ld_{1}+\ld_2},\ldots, s_{\ld_{1} + \ldots + \ld_{n-1}} \}.$
For each $z\in \ZZ$, let $\ld + z$ be the composition in $\Ld$ such that
$\W_{\ld+z} = \pi^{-z}\W_\ld \pi^z$.
\exa
Let $n=3, d=6, z=4$ and $\ld = (1,2,3)\in\Ld$. We have $W_{\ld} = \g{s_2, s_4, s_5, s_6}$, $W_{\ld+4} = \g{s_6, s_2,s_3,s_4}$ and hence $\ld+4 = (1,4,1)\in\Ld$.
\endexa
Let $\D_\ld = \{w\in \tW ~|~ \ell(wg) = \ell(w) + \ell(g) \fa g\in \W_\ld\}$.
Then $\D_\ld$ (resp., $\D_\ld^{-1}$) is the set of distinguished right (resp. left) coset representatives of $\W_\ld$ in $\tW$.
Denote by $\D_{\ld\mu} = \D_\ld \cap \D_\mu^{-1}$ the set of distinguished double coset representatives.
\begin{lem}\label{lem:doublecoset}
Let $\ld,\mu \in \Ld$, and let $g \in \D_{\ld\mu}$. Then
\enua
\item There is a unique $\delta =\delta(\ld,g,\mu) \in \Ld(n',d)$ for some $n'$ such that
\[
W_{\delta} = g^{-1} W_\ld g \cap W_\mu.
\]
\item The map $W_\ld \times (\D_\delta \cap W_\mu) \rw W_\ld g W_\mu$ sending $(x,y)$ to $xgy$ is a bijection satisfying $\ell(xgy) = \ell(x) + \ell(g) + \ell(y)$.
\endenua
\end{lem}
\proof
Part (a) follows from \cite[Lemma~2.2.2]{Gr99} and Part (b) is known (cf. \cite[Theorem~4.18]{DDPW08}).
\endproof
Let $\leq$ be the (strong) Bruhat order on $\W$. Extend it to $\tW$ by
\[
\pi^{z_1} w_1 \leq \pi^{z_2} w_2
\textup{ if and only if }
z_1 = z_2, w_1 \leq w_2.
\]
\subsection{Extended affine Hecke algebras}
The extended affine Hecke algebra $\cH = \cH(\tW)$ associated to $\tW$ is a $\ZZ[v,v^{-1}]$-algebra with a basis $\{T_g ~|~ g\in\tW\}$
(cf. e.g., \cite[Proposition~1.2.3]{Gr99}) satisfying $T_w T_{w'} = T_{ww'}$ if $\ell(w) + \ell(w') = \ell(ww')$ and $(T_{s} +1)(T_{s}  - v^2) = 0$ for $s\in S$.
For a finite subset $X \subset \tW$ and for each $\ld\in\Ld$, let
\[
T_X = \sum_{w\in X} T_w
\quad
\textup{and}
\quad
x_\ld = T_{\W_\ld}.
\]
Following \cite{KL79}, denote by $\{C'_w~|~ w\in \W\}$ the Kazhdan-Lusztig basis of the Hecke algebra $\cH(\W)$ associated to $W$. For each $w\in W$, we have
$C'_w = v^{-\ell(w)} \sum_{y \leq w} P_{y,w}T_y$,
where $P_{y,w} \in \ZZ[v^2]$ is the Kazhdan-Lusztig polynomial.
Note that $\cH = \cH(\tW)$ contains $\cH(\W)$ as a subalgebra, we define $C'_g = T_\pi^z C'_w \in \cH$ for each $g = \pi^z w \in \tW$ with $w\in W, z\in\ZZ$.

Statements in Lemma~\ref{lem:Cur} below are known for non-extended Weyl groups and Hecke algebras (cf. \cite[Theorem~1.2(i)]{Cur85}, \cite[Corollary 4.19]{DDPW08}). It seems that the extended version is taken for granted for the experts. Statements (a) and (c) can be found in \cite[(7.1.1)]{DF14} and \cite[Lemma~3.1]{FS14}, respectively.
\begin{lemma}\label{lem:Cur}
Let $\ld,\mu \in \Ld, g = \pi^z w\in \D_{\ld\mu}$ for some $w\in \W$ and $z\in \ZZ$.
Denote by $w_\circ^\nu$ the longest element in $\W_\nu$ for any composition $\nu$. Then:
\enua
\item The longest element $g_{\ld\mu}^+$ in $W_\ld g W_\mu$ is given by $g_{\ld\mu}^+ = w_\circ^\ld g w_\circ^{\delta(\ld,g,\mu)} w_\circ^\mu$. In particular,
\[
\ell(g_{\ld\mu}^+) =  \ell(w_\circ^\ld) +  \ell(g)  - \ell(w_\circ^{\delta(\ld,g,\mu)}) +  \ell(w_\circ^\mu).
\]
\item $W_\ld g W_\mu = \{x \in \tW ~|~ g \leq x \leq g^+_{\ld\mu}\}$.
\item There exists $c_{x,g}^{(\ld,\mu)}\in\ZZ[v,v^{-1}]$ such that
\[
T_{W_\ld g W_\mu}
= v^{\ell(g^+_{\ld\mu})} C'_{g^+_{\ld\mu}} + \sum_{\substack{x\in \D_{\ld\mu} \\x < g}} c_{x,g}^{(\ld,\mu)} C'_{x^+_{\ld\mu}}.
\]
In particular,
$
x_\mu
=
v^{\ell(w_\circ^\mu)}
C'_{w_\circ^\mu}.$
\endenua
\end{lemma}

\subsection{Affine $q$-Schur algebras}\label{sec:S}

For $\ld,\mu\in\Ld$ and $g\in \D_{\ld\mu}$,
denote by $\phi_{\ld\mu}^g\in \Hom_\cH(x_\mu\cH, x_\ld\cH)$ the right $\cH$-linear map sending  $x_\mu$ to $T_{W_\ld g W_\mu}$.
Thanks to Lemma~\ref{lem:doublecoset}(b),
we have
$T_{W_\ld g W_\mu} = x_\ld T_g T_{\D_\delta \cap W_\mu}$ for some
$\delta \in \Ld(n',d)$ and hence $\phi_{\ld\mu}^g \in \Hom_\cH(x_\mu\cH,
x_\ld\cH)$. The affine $q$-Schur algebra is defined by
\[
\Sq =\Sq(n,d) = \textup{End}_{\cH}
\Bp{
\mathop{\oplus}_{\ld\in\Ld} x_\ld \cH
}
= \bigoplus_{\ld,\mu \in \Ld} \textup{Hom}_{\cH} (x_\mu \cH, x_\ld \cH)
.
\]
There is also a geometric definition for $\Sq$ as given in \cite{Lu99}.
It is known (cf. \cite[Theorem~2.2.4]{Gr99}) that $\{  \phi_{\ld\mu}^g ~|~ \ld,\mu \in \Ld, g \in \D_{\ld\mu} \}$ forms a basis of $\Sq$.

Let $\Tt = \bigcup_{d \in \NN} \Tt_d$, where $\Tt_d$ is the set of $\ZZ\times\ZZ$ matrices over $\NN$ in which each element $A = (a_{ij})_{ij}$ satisfies the following conditions:
\enu
\item[(T1)] $a_{ij} = a_{i+n,j+n}$ for all $i,j \in \ZZ$;
\item[(T2)] $\sum\limits_{i=1}^{n}\sum\limits_{j\in\ZZ} a_{ij} = d$.
\endenu
For $i,j\in\ZZ$, define a matrix  $E_{ij} \in \Tt_1$ by
\[
(E_{ij})_{xy} =\bc{
1&\tif (x,y) = (i+rn,  j+rn) \textup{ for the same } r\in\ZZ,
\\
0&\otw.
}
\]
 For each matrix $T = (t_{ij})_{ij}\in \Tt$, define its  row sum vector
 $\ro(T)= (\ro(T)_1, \ldots, \ro(T)_n)$
and its column sum vector $\co(T)= (\co(T)_1, \ldots \co(T)_n)$ by
\[
\ro(T)_k = \sum\limits_{j\in\ZZ} t_{kj},
\quad
\co(T)_k = \sum\limits_{i\in\ZZ} t_{ik},
\quad
k = 1, \ldots, n.
\]
For each $\ld \in \Ld$ and $i = 1,\ldots, n$, we define integral intervals with respect to $\lambda$ by
\[
R_i^\ld =
(\textstyle\sum_{k=1}^{i-1}\ld_k ~..~ \sum_{k=1}^{i}\ld_k].
\]
The basis for $\Sq$ can be parametrized by $\Tt_d$ thanks to the following identities (cf. \cite[\S 7.4]{VV99}, \cite{DF15}).
\begin{lemma}\label{lem:kappa}
The map
\[
\kappa: \{(\ld,g,\mu) ~|~ \ld,\mu\in\Ld(n,d), g\in\D_{\ld\mu}\} \longrightarrow \Tt_d
\]
is a bijection sending $(\ld,g,\mu)$ to $A = (a_{ij})_{ij}$ where $a_{ij} = |R_i^\ld \cap gR_j^\mu|$.
\end{lemma}
For $A = (a_{ij})_{ij} = \kappa(\ld,g,\mu) \in \Tt_d$, set $e_A = \phi_{\ld\mu}^g$.
Hence $\{e_A ~|~ A\in \Tt_d\}$ forms a basis of $\Sq$.
For each $j = 1, \ldots, n$, let $(\delta^{(j)}_{1}, \ldots, \delta^{(j)}_{k_j}) \in \Ld(k_j, \ld_j)$ for some $k_j \in \NN$ be the composition obtained from $(\ldots, a_{-1,j}, a_{0j}, a_{1j},  \ldots)$ by deleting all zero entries.
Define $\delta(A) \in \Ld(n',d)$ by
\eq\label{eq:delta}
\delta(A) = (\delta^{(1)}_{1}, \ldots, \delta^{(1)}_{k_1}, \delta^{(2)}_1, \ldots, \delta^{(n)}_1, \ldots, \delta^{(n)}_{k_n}).
\endeq
The Following lemma can be found in \cite[Corollary 3.2.3]{DDF12}.
\begin{lemma}
Let $A = \kappa(\ld,g,\mu)$. Then
$W_{\delta(A)} = g^{-1} W_\ld g \cap
W_\mu$.
In particular, $\delta(A)$ is equal to $\delta(\ld,g,\mu)$ described
in Lemma~\ref{lem:doublecoset}.
\end{lemma}

Set $\ell(A) = \ell(g)$ for $A = \kappa(\ld,g,\mu)$. The following lemma and a proof can be found in \cite[Lemma~3.2(2)]{DF15}. Their proof also provided the minimal representative of a double coset associated to $A$. Here we provide a proof by combining Lemma~\ref{lem:l(g)} and Lemma~\ref{lem:kappa}.
\begin{lem}\label{lem:l(A)}
Assume that $A\in \Tt_d$. Then
\[
\ell(A)
=
\sum\limits_{\substack{i \in \ZZ \\ 1\leq j \leq n} }
\sum\limits_{\substack{x < i\\ y > j}}
a_{ij} a_{xy}
=
\sum\limits_{\substack{1\leq i \leq n \\ j \in \ZZ} }
\sum\limits_{\substack{x >i \\ y < j}}
a_{ij} a_{xy}.
\]
\end{lem}
\proof
Let $A = \kappa(\ld,g,\mu)$ for some $\ld,\mu\in \Ld$ and $g\in \D_{\ld\mu}$.
By Lemma~\ref{lem:kappa}, for all $i,j \in \ZZ$ there is a natural bijection
$R_i^\ld \cap gR_j^\mu \leftrightarrow \{(g(s),s) \in R_i^\ld\times R_j^\mu\}$ between sets of size $a_{ij}$. Note that for $(g(s),s), (g(t),t) \in R_i^\ld\times R_j^\mu$, the condition $s<t$ is equivalent to the condition $g(s) < g(t)$ since $g\in \D_{\ld\mu}$. Hence if $(s,t) \in R^\mu_j \times R^\mu_y$ satisfies both $s < t$ and $g(s) > g(t)$, then $j$ must be smaller than $y$.

Under these bijections, the set of pairs $(s,t)\in \ZZ^2$ satisfying ``$s<t, g(s)>g(t)$ and $s\in [1..d]$'' becomes the set of quadruples $(g(s),s, g(t),t) \in R_i^\ld \times R_j^\mu \times R_x^\ld \times R_y^\mu$ satisfying ``$j<y, i>x$ and $1\leq j \leq n$''.
 The first assertion follows. The second assertion follows from that $\ell(g) = \ell(g^{-1})$ and $\kappa(\mu,g^{-1},\ld) = \ltrans{A}$.
\endproof
\exa
Let $n=2, d=10$ and
\[
A = 3E_{10} + 4E_{12} + E_{23} + 2E_{24} =
\scalebox{0.85}{$
\bM{[cc:cc:cccc]
\ddots&4&&&&\\
&&1&2&\\
\hdashline
&3&0&4&\\
&&0&0&1&2&\\
\hdashline
&&&3&&4\\
&&&&&\ddots\\
}
$}
 \in \Tt_{10}.
\]
We have $\delta(A) = (a_{01}, a_{02},a_{12},a_{32}) = (1,2,4,3) \in \Ld(4,10)$ and
$\ell(A) = 3(1+2) = 9$.
\endexa
\subsection{The standard basis}
For $A \in \Tt_d$, let
\[
d_A =
\sum_{\substack{1\leq i \leq n\\j \in \ZZ} }
\sum_{\substack{x \leq i\\ y > j}}
a_{ij} a_{xy}
,\qquad
[A] = v^{-d_A} e_A.
\]
It is clear that $\{[A] ~|~ A \in \Tt_d\}$ is a basis of $\Sq$,
which is called the \textit{standard basis} (cf. \cite{Lu99}).
The following lemma is due to Du-Fu \cite[Lemma~7.1]{DF14},
and we offer a slightly different argument.
\lem\label{lem:dA}
For $A =\kappa(\ld,g,\mu)\in \Tt_d$, we have
$
d_A = \ell(g^+_{\ld\mu})-\ell(w_\circ^\mu).
$
\endlem
\proof
Let $\delta = \delta(\ltrans{A}) = (\delta^{(1)}_{1}, \ldots, \delta^{(1)}_{k_1}, \delta^{(2)}_1, \ldots, \delta^{(n)}_1, \ldots, \delta^{(n)}_{k_n})$ as in \eqref{eq:delta}.
So $\ld_i = \sum_{j=1}^{k_i} \delta_j^{(i)}$,
$\W_\delta \simeq \W_{\delta(A)}$ and hence
$\ell(\wo^{\delta(A)}) = \ell(\wo^{\delta})$.
We have
$\ell(g_{\ld\mu}^+) - \ell(w_\circ^\mu) =  \ell(g) + \ell(w_\circ^\ld) - \ell(w_\circ^{\delta})$
by Lemma~\ref{lem:Cur}(a), where
\[
\ba{{rlll}
\ell(g) &\ds= \sum_{\substack{1\leq i \leq n\\j \in \ZZ} }
\sum_{\substack{x < i\\ y > j}}
a_{ij} a_{xy},
\\
\ell(w_\circ^\ld) - \ell(w_\circ^{\delta})
\ds
&\ds=
 \sum_{i=1}^n \tbinom{\ld_i}{2}
- \sum_{i=1}^{n'} \tbinom{\delta_i}{2}
=
\sum_{i=1}^n
\Bp{
\tbinom{\sum_{j=1}^{k_i} \delta_{j}^{(i)}}{2}
-  \sum_{j=1}^{k_i} \tbinom{\delta_{j}^{(i)}}{2}
}
=\ds \sum_{i=1}^n \sum_{y>j} \delta_{j}^{(i)} \delta_{y}^{(i)}
\\
&\ds=
\sum_{\substack{1\leq i \leq n\\j \in \ZZ} }
\sum_{\substack{x = i\\ y > j}}
a_{ij} a_{xy}.
}
\]
Therefore,
$
\ell(g_{\ld\mu}^+) - \ell(w_\circ^\mu) =
\sum\limits_{\substack{1\leq i \leq n\\j \in \ZZ} }
\sum\limits_{\substack{x \leq i\\ y > j}}
a_{ij} a_{xy} = d_A.
$
\endproof
Denote the bar involution on $\cH$ by
$
\bar{~}:\cH \rw \cH,\quad v \mapsto v^{-1},\quad T_g \mapsto T_{g^{-1}}^{-1}.
$
Following \cite[Proposition~3.2]{Du92}, the bar involution on $\Sq$ can be described as follows:
for each $f \in \textup{Hom}_{\cH}(x_\mu \cH, x_\ld \cH)$, let $\={f}\in\textup{Hom}_{\cH}(x_\mu \cH, x_\ld \cH)$ be the map sending $v$ to $v^{-1}$ and $C'_{w_\circ^\mu}$ to $\={f(C'_{w_\circ^\mu})}$.
Equivalently,
\[
\={f}(x_\mu H) = v^{2\ell(w_\circ^\mu)}\={f(x_\mu)} H ~\fa~ H \in \cH.
\]
In particular, for $A=\kappa(\ld,g,\mu) \in \Tt_d$, by Lemma~\ref{lem:Cur} we have
\eqnarray
e_A(C'_{w_\circ^\mu})
&\ds = v^{\ell(g^+_{\ld\mu})-\ell(w_\circ^\mu)} C'_{g^+_{\ld\mu}}
 + \sum_{\substack{x\in\D_{\ld\mu}\\x < g}}
v^{-\ell(w_\circ^\mu)} c_{x,g}^{(\ld,\mu)} C'_{x^+_{\ld\mu}},
    \label{eq:eA}
\\
\={e_A}(C'_{w_\circ^\mu})
&\ds = v^{\ell(w_\circ^\mu)-\ell(g^+_{\ld\mu})} C'_{g^+_{\ld\mu}}
+ \sum_{\substack{x\in\D_{\ld\mu}\\x < g}}
v^{\ell(w_\circ^\mu)}\={c_{x,g}^{(\ld,\mu)}} C'_{x^+_{\ld\mu}}.
    \label{eq:eAbar}
\endeqnarray
\prop\label{prop:Abar}
Assume that $A =\kappa(\ld,g,\mu)\in \Tt_d$. There exists $\gamma^{(\ld,\mu)}_{x,g}  \in \ZZ[v,v^{-1}]$ for each $x\in \D_{\ld\mu}$ such that
\[
\={[A]} = [A] + \sum_{\substack{x\in \D_{\ld\mu} \\ x< g}} \gamma^{(\ld,\mu)}_{x,g} [\kappa(\ld,x,\mu)].
\]
\endprop
\proof
By Lemma~\ref{lem:dA}, Equations \eqref{eq:eA} and \eqref{eq:eAbar} can be rewritten as
\begin{eqnarray*}
[A](C'_{w_\circ^\mu})
&\ds =  C'_{g^+_{\ld\mu}}
 + \sum_{\substack{x\in\D_{\ld\mu}\\x < g}}
v^{-\ell(g^+_{\ld\mu})} c_{x,g}^{(\ld,\mu)} C'_{x^+_{\ld\mu}},
\\
\={[A]}(C'_{w_\circ^\mu})
&\ds =  C'_{g^+_{\ld\mu}}
+ \sum_{\substack{x\in\D_{\ld\mu}\\x < g}}
v^{\ell(g^+_{\ld\mu})}\={c_{x,g}^{(\ld,\mu)}} C'_{x^+_{\ld\mu}}.
\end{eqnarray*}
If $\ell(g) = 0$ (i.e. $g= \pi^z$ for some $z$) then $\={[A]} = [A]$
and we are done. For arbitrary $g$, it follows from an easy
induction on $\ell(g)$.
\endproof
A matrix $A = (a_{ij})_{ij}$ is called \textit{bidiagonal} if either $a_{ij} = 0$ for all $j \neq i, i+1$ or  $a_{ij} = 0$ for all $j \neq i, i-1$.
\begin{cor}\label{cor:Abi}
\textup{(}\cite[Lemma~7.2]{DF14}\textup{)} If $A \in \Tt_d$ is
bidiagonal then $[A]$ is bar-invariant.
\end{cor}
\proof
By Lemma~\ref{lem:l(A)}, $\ell(A) = 0$ for any bidiagonal matrix $A$ and we are done.
\endproof
We define a partial order $\leq_a$ on $\Tt$ by $A \leq_a B$ if and only if $\ro(A) = \ro(B)$, $\co(A)=\co(B)$ and $\sigma_{i,j}(A) \leq \sigma_{i,j}(B)$ for all $i \neq j$ where
\[
\sigma_{i,j}(A) =
\bc{
\sum\limits_{x\leq i, y\geq j} a_{xy} &\tif i > j,
\\
\sum\limits_{x\geq i, y\leq j} a_{xy} &\tif i < j.
}
\]
In the following the expression ``lower terms'' represents a linear combination of smaller elements with respect to the partial order $\leq_a$. Here is an affine generalization of \cite[Lemma~3.6]{BLM90}.
\begin{lemma}\label{lem:PO}
\cite[Lemma 3.6]{FS14}
Assume that $A = \kappa(\ld, g, \mu)$ and $B = \kappa(\ld, h, \mu)$. If $h \leq g$  then $B \leq_a A$.
\end{lemma}

Below is a corollary obtained by combining Proposition~\ref{prop:Abar} and Lemma~\ref{lem:PO}. This statement has also appeared in the proof of \cite[Proposition 7.6]{DF14}.
\begin{cor}\label{cor:Abar}
For $A \in \Tt_d$, we have
\[
\={[A]} = [A] + \textup{lower terms}.
\]
\end{cor}

\section{Construction of a monomial basis}
\subsection{Multiplication formulas}
For each $A \in \Tt$, let $\textup{diag}(A) = (\delta_{ij} a_{ij})_{ij} \in \Tt$ and let $A^\pm \in \Tt$ be such that
\[
A = A^\pm + \textup{diag}(A).
\]
For any matrix $T= (t_{ij})_{ij} \in \Tt$, denote the matrix obtained by shifting every entry of $T$ up by one row as
\[
\^{T}= (\^{t}_{ij})_{ij}, \quad \^{t}_{ij} = t_{i+1,j}.
\]
On the other hand, denote the matrix obtained by shifting every entry of $T$ down by one row as
\[
\V{T}= (\V{t}_{ij})_{ij}, \quad \V{t}_{ij} = t_{i-1,j}.
\]
For $A= (a_{ij})_{ij}, B=(b_{ij})_{ij} \in \Tt$, define
\[
\ds\lrb{A+B}{A} =
\prod_{\substack{1\leq i \leq n \\ j\in \ZZ}}
\dfrac{ [a_{ij}+b_{ij}][a_{ij}+b_{ij}-1] \ldots [b_{ij}+1] }
{ [a_{ij}][a_{ij}-1] \ldots [1] }.
\]
The following remarkable multiplication formulas were due to \cite[Proposition~3.6]{DF15}.
\begin{lem}\label{lem:mult}
Assume that $A,B \in \Tt_d$,
$\ro(A) = \co(B)$ and $B$ is bidiagonal.
Let $\Tt_\af = \{ T \in \Tt ~|~ \ro(T) = \af \}$ for $\af \in \Ld$.
\enua
\item If $B$ is upper triangular (i.e. $B^\pm = \sum \af_i E_{i-1,i}$), then
\eq\label{eq:multA}
[B]* [A] = \sum_{T \in \Tt_\af}
v^{\beta(A,T)}
\={\lrb{A - T + \^{T}}{A-T} }
[A - T + \^{T}],
\endeq
where
\[
    \beta(A,T) =
    \sum\limits_{1 \leq i\leq n } \sum\limits_{ j\leq y }
    \^{t}_{ij} (a_{iy} - t_{iy})
    -
    \sum\limits_{1 \leq i\leq n }\sum\limits_{j<y}
    t_{ij}(a_{iy}-t_{iy}).
\]
\item If $B$ is lower triangular (i.e. $B^\pm = \sum \af_i E_{i+1,i}$), then
\eq\label{eq:multB}
[B]* [A] = \sum_{T \in \Tt_\af}
v^{\beta'(A,T)}
\={\lrb{A  - T + \V{T}}{A-T} }
[A  - T + \V{T}],
\endeq
where
\[
    \beta'(A,T) =
    \sum\limits_{1 \leq i\leq n }\sum\limits_{ j\geq y }
    \V{t}_{ij} (a_{iy} - t_{iy})
    -
    \sum\limits_{1 \leq i\leq n }\sum\limits_{j>y}
    t_{ij}(a_{iy}-t_{iy}).
\]
\endenua
\end{lem}
\alg\label{alg:ht} Assume that $A,B \in \Tt_d$,
$\ro(A) = \co(B)$ and $B$ is bidiagonal. We produce a matrix $M
\in \Tt_d$ as follows. \enua
\item If $B$ is upper triangular (i.e. $B^\pm = \sum \af_i E_{i-1,i}$), then:
\enu
\item[(1)] For each row $i$, find the unique $j$ such that $\af_i \in \big{(} \sum_{y>j} a_{iy} .. \sum_{y\geq j} a_{iy} \big{]}$.
\item[(2)] Construct a matrix
$T_+ = \sum_{i=1}^n \bp{ (\af_i - \sum_{y>j} a_{iy}) E_{ij} + \sum_{y>j} a_{iy} E_{iy}}$.
\item[(3)] Let $M = A - T_+ + \^{T}_+$.
\endenu
\item If $B$ is lower triangular (i.e. $B^\pm = \sum \af_i E_{i+1,i}$), then:
\enu
\item[(1)] For each row $i$, find the unique $j$ such that $\af_i \in \big{(} \sum_{y<j} a_{iy} .. \sum_{y\leq j} a_{iy} \big{]}$.
\item[(2)] Construct a matrix $T_+ = \sum_{i=1}^n \bp{ (\af_i - \sum_{y<j} a_{iy} )E_{ij} + \sum_{y<j} a_{iy} E_{iy} }$.

\item[(3)] Let $M = A - T_+ + \V{T}_+$.
\endenu
\endenua
That is, the matrix $M$ is obtained from $A$ by  ``shifting'' up (or down) entries by one row
starting from the rightmost (or leftmost) nonzero entries on each row.
\endalg
\lem\label{lem:ht}
The highest term (with respect to $\leq_a$) in \eqref{eq:multA} or in \eqref{eq:multB} exists and its corresponding matrix is the matrix $M$ described in Algorithm \ref{alg:ht}.
\endlem
\proof
If $B$ is upper triangular, then each term on the right-hand side of \eqref{eq:multA} must be of the form $[A-T+\^{T}]$ for some $T\in \Tt_\af$ such that $a_{ij}-t_{ij}+\^{t}_{ij} \geq 0$ for all $i,j \in \ZZ$.
Note that
\[
\sigma_{ij}(A-E_{xy}+\^{E}_{xy})
=\bc{
\sigma_{ij}(A)+1&\tif j< i=x-1, j\leq y,
\\
\sigma_{ij}(A)-1&\tif j>i=x, j\geq y,
\\
\sigma_{ij}(A)&\otw.
}
\]
It follows immediately that, for each $i$,
\[
\ldots <_a (A-E_{i,-1}+\^{E}_{i,-1}) <_a (A-E_{i0}+\^{E}_{i0}) <_a (A-E_{i1}+\^{E}_{i1}) <_a \ldots
\]
Therefore, for any $T\in\Tt_\af$ we have $A-T+\^{T} \leq_a A-T_++\^{T}_+ = M$.

The case that $B$ is lower triangular is similar and skipped.
\endproof
\exa
Let $n = 2$, $B^\pm = 2E_{12} + 1E_{23}$ and $A = 2E_{12}+3E_{21}+E_{22}+E_{23}$, that is,
\[
B = \bM{[c:cc:cc]
\ddots&1&&\\
\hdashline
&*&2&\\
&&*&1\\
\hdashline
&&&\ddots
},
\quad
A = \bM{[c:cc:cc]
\ddots&1&&\\
\hdashline
&&2&\\
&3&1&1\\
\hdashline
&&&\ddots
}.
\]
Then $\af_1 = 1 \in \Big{(} \sum\limits_{y>2} a_{1y} .. \sum\limits_{y\geq 2} a_{1y} \Big{]} = (0..2]$ and
$\af_2 = 2 \in \Big{(} \sum\limits_{y>2} a_{2y} .. \sum\limits_{y\geq 2} a_{2y} \Big{]} = (1..2]$. Therefore
\[
T_+= \bM{[c:cc:cc]
\ddots&1&&\\
\hdashline
&&1&\\
&&1&1\\
\hdashline
&&&\ddots
},
\quad
M = \bM{[c:cc:cc]
\ddots&\uca{0}&1&\\
\hdashline
&&\uca{2}&1\\
&3&\uca{0}&\uca{0}\\
\hdashline
&&&\ddots
}.
\]
\endexa
\subsection{Admissible pairs}
We call a pair $(B,A)$ of matrices to be \textit{admissible} if either of the following conditions (A1) or (A2) holds.
\enu
\item[(A1)]
$
\ba{[t]{lll}
B^\pm &= \sum\limits_{i=1}^n m_{i} E_{i,i+1}\fsm m_i \in \NN,\textup{ and}
\\
A^\pm &= \sum\limits_{i = 1}^n \sum\limits_{j\leq k} a_{i,i+j} E_{i,i+j} \fsm k \in \ZZ,
\textup{ where }a_{i,i+k} \geq m_i \fa i;
}
$
\item[(A2)]
$((b_{-i,-j})_{ij}, (a_{-i,-j})_{ij})$ satisfies Condition (A1).
\endenu
That is, if $(B,A)$ satisfies Condition (A1), we have
\[
\scalebox{0.85}{$
B =\bM{[c:cccc]
\ddots&&\\
\ddots&m_{n}\\
\hdashline
&*&m_{1}\\
&&*&m_2&\\
&&&*&\ddots\\
&&&&\ddots
},~
A =\bM{[ccccc]
&&\\
\ddots&\\
\hdashline
\ddots&*+m_n&\\
\ddots&*&*+m_1&\\
\ddots&*&*&\ddots\\
\ddots&\ddots&\ddots&\ddots
},~
M =\bM{[ccccc]
\ddots&&\\
\ddots&m_n&\\
\hdashline
\ddots&\uca{*}&m_1&\\
\ddots&*&\uca{*}&\ddots\\
\ddots&*&*&\ddots\\
\ddots&\ddots&\ddots&\ddots
}.$}
\]
\thm\label{thm:M}
If $(B,A)$ is admissible then
$
[B] * [A]
=
[M]+ \textup{lower terms}.
$
\endthm
\proof
We only prove when $B$ is upper triangular since the other case is similar.
Due to Lemma~\ref{lem:ht}, it remains to show that the coefficient for $[M]$ is one.
If $(B,A)$ is admissible, then $T_+ =  \sum\limits_{i=1}^n m_i E_{i,i+k}$ and hence
\[
\lrb{A-T_++\^{T}_+}{A-T_+}
=
\prod_{1\leq i \leq n}
\left(
\prod_{j < i+k}
\dfrac{[(A-T_+)_{ij}] \ldots [1]}{[(A-T_+)_{ij}] \ldots [1]}
\right)
= 1.
\]
Note that by definition of admissible pairs we have $\sum_{ j\leq y } (a_{iy} - t_{iy}) = 0$ for each nonzero $\^{t}_{ij}$ and
$\sum_{j<y} (a_{iy}-t_{iy}) = 0$ for each nonzero $t_{ij}$. Hence $\beta(A,T_+) = 0$.
\endproof
\subsection{A monomial basis}
Below we provide an algorithm that generates a monomial basis in a diagonal-by-diagonal manner involving only  admissible pairs (see also \cite{FL14} for a diagonal-by-diagonal construction in a finite type setting).
\alg\label{alg:mono}
For each $A=(a_{ij})_{ij} \in \Tt_d$,
we construct upper bidiagonal matrices $B^{(1)}, \ldots, B^{(x)}$ and lower bidiagonal matrices $B_{(1)},\ldots, B_{(y)}$ as follows:
\enu
\item Initialization: $t=0$, $U^{(0)} = A$.
\item If $U^{(t)}$ is a lower triangular matrix, then go to Step (5) (denote this $t$ by $x$). Otherwise, denote the outermost nonzero upper diagonal of the matrix $U^{(t)}= (u^{(t)}_{ij})_{ij}$ by
$T_+^{(t)}=   \sum_{i=1}^n  u_{i,i+k}^{(t)} E_{i,i+k}$ for some $k>0$.
\item Define matrices
\[
\ba{{ll}
\ds B^{(t+1)} = \sum_{i=1}^n u^{(t)}_{i,i+k} E_{i,i+1} + \textup{a diagonal determined by }\eqref{eq:diag},
\\
U^{(t+1)} = U^{(t)} - T_+^{(t)} + \V{T}_+^{(t)}.
}
\]
\item Increase $t$ by one and then go to Step (2).
\item Set $L^{(0)} = U^{(x)}$ and set $s=0$.
\item If $L^{(s)}$ is a lower bidiagonal matrix (denote this $s$ by $y$), then set $B_{(y)} = L^{(y)}$ and end the algorithm.
Otherwise, denote the outermost nonzero lower diagonal of the matrix $L^{(s)} = (l^{(s)}_{ij})_{ij}$ by $T_{+,(s)}=   \sum_{i=1}^n  l^{(s)}_{i+k,i} E_{i+k,i}$ for some $k>0$.
\item Define matrices
\[
\ba{{lll}
\ds B_{(s+1)} = \sum_{i=1}^n l^{(s)}_{i+k,i} E_{i+1,i} + \textup{a diagonal determined by }\eqref{eq:diag},
\\
\ds L^{(s+1)} = L^{(s)} - T_{+,(s)} + \^{T}_{+,(s)}.
}
\]
\item Increase $s$ by one and then go back to Step (6).
\endenu
Here the diagonal entries are uniquely determined by
\eq\label{eq:diag}
\ba{{lll}
\ro(B^{(1)}) = \ro(A),
&\co(B^{(i)}) = \ro(B^{(i+1)})
&\textup{ for } i = 1,\ldots, x-1,
\\
\co(B^{(x)}) = \ro(B_{(1)}),
&\co(B_{(i)}) = \ro(B_{(i+1)})
&\textup{ for } i = 1,\ldots, y-1.
}
\endeq
\endalg
\thm\label{thm:min gen}
For $A \in \Tt_d$, the matrices $B^{(1)}, \ldots, B^{(x)}, B_{(1)}, \ldots, B_{(y)} \in \Tt_d$ in Algorithm~\ref{alg:mono} satisfy that
\[
[ B^{(1)}]
*
\cdots
*
[ B^{(x)}]
*
[B_{(1)}]
*
\cdots
*
[B_{(y)}]
= [A] +\textup{lower terms}.
\]
\endthm
\proof
For each admissible pair $(Y,X)$, let $M$ be the matrix corresponding to the highest term in $[Y]*[X]$ (cf. Algorithm~\ref{alg:ht}). For any matrix  $X' <_a X$, let $M'$ be the matrix that corresponds to the highest term in $[Y]*[X']$. By construction we have $M' <_a M$, and hence
\[
[Y]*([X] + \textup{lower terms}) = [M] + \textup{lower terms}.
\]
Algorithm \ref{alg:mono} guarantees that each pair $(B^{(j)}, U^{(j)})$ or $(B_{(j)}, L^{(j)})$ is admissible (here it is understood that $L^{(y-1)} = B_{(y)}$ and $U^{(x-1)} = L^{(0)}$). Hence by Theorem~\ref{thm:M},
\[
\ba{{ll}
[ B^{(1)}] * \cdots * [ B^{(x-1)}] * [ B_{(1)}] * [ B_{(2)}] * \cdots * \bp{ [ B_{(y-1)}] * [ B_{(y)}] }
\\
= [ B^{(1)}] * \cdots * [ B^{(x-1)}] * [ B_{(1)}] * [ B_{(2)}] * \cdots * [ B_{(y-2)}] * \bp{ [ L^{(y-2)}] + \textup{lower terms}}
\\
= \ldots
\\
=
[ B^{(1)}] * \cdots * [ B^{(x-1)}] * \bp{ [L^{(0)}] + \textup{lower terms}}
\\
=
[ B^{(1)}] * \cdots * [ B^{(x-1)}] * \bp{ [U^{(x-1)}] + \textup{lower terms}}
\\
= \ldots
\\
= [A] + \textup{lower terms}.
}
\]
\endproof
For each $A\in \Tt_d$, we define
\eq\label{eq:mono}
m_A = [ B^{(1)}]
*
\cdots
*
[ B^{(x)}]
*
[B_{(1)}]
*
\cdots
*
[B_{(y)}]
.
\endeq
\begin{cor}\label{cor:mA}
The set
$\{m_A ~|~ A \in \Tt_d\}$ forms a basis of the $\ZZ[v,v^{-1}]$-algebra $\Sq$ (called a \textit{monomial basis}).
Moreover, $m_A$ is bar invariant for each $A\in\Tt_d$.
\end{cor}
\proof
The first assertion is clear from Theorem~\ref{thm:min gen}. The second assertion follows from Corollary \ref{cor:Abi}.
\endproof
\rmk
In \cite{DDX07}, Deng, Du and Xiao constructed a family of monomial bases for the Hall algebra of the cyclic quiver. Any such basis can be adapted to a monomial basis for $\Sq$ using surjections from the double Hall algebra of the cyclic quiver to $\Sq$ (cf. \cite{DF14}).
But the relation between their monomial bases and ours is unclear.
\endrmk
\subsection{An example}
\exa
Let $n = 2, d=21$, and let
\[
A =E_{01}+2E_{02}+3E_{03}+4E_{12}+5E_{21}+6E_{20} = \bM{[cc:cc:cc]
\ddots&&&&&\\
&0&1&2&3\\
\hdashline
&0&0&4&0\\
&6&5&0&1\\
\hdashline
&0&0&0&0\\
&&&&&\ddots
}\in \Tt_{21}.
\]
We have $[A] = [B^{(1)}]*[U^{(1)}] + $ lower terms, where
\[
B^{(1)} =   \bM{[cc:cc:cc]
\ddots&\ddots&&&&\\
&*&3&&\\
\hdashline
&&*&0&\\
&&&*&3\\
\hdashline
&&&&*&\ddots\\
&&&&&\ddots
},
\quad
U^{(1)} =   \bM{[cc:cc:cc]
\ddots&&&&&\\
&0&1&2&\\
\hdashline
&0&0&4&\dca{3}\\
&6&5&0&1\\
\hdashline
&0&0&0&0\\
&&&&&\ddots
}.
\]
$[U^{(1)}] = [B^{(2)}]*[U^{(2)}] + $ lower terms, where
\[
B^{(2)} =   \bM{[cc:cc:cc]
\ddots&\ddots&&&&\\
&*&2&&\\
\hdashline
&&*&3&\\
&&&*&2\\
\hdashline
&&&&*&\ddots\\
&&&&&\ddots
},
\quad
U^{(2)} =   \bM{[cc:cc:cc]
\ddots&&&&&\\
&0&\dca{4}&&\\
\hdashline
&0&0&\dca{6}&\\
&6&5&0&\dca{4}\\
\hdashline
&0&0&0&0\\
&&&&&\ddots
}.
\]
$[U^{(2)}] = [B^{(3)}]*[U^{(3)}] + $ lower terms, where
\[
B^{(3)} =   \bM{[cc:cc:cc]
\ddots&\ddots&&&&\\
&*&4&&\\
\hdashline
&&*&6&\\
&&&*&4\\
\hdashline
&&&&*&\ddots\\
&&&&&\ddots
},
\quad
U^{(3)} =   \bM{[cc:cc:cc]
\ddots&&&&&\\
&6&&&\\
\hdashline
&0&\dca{4}&&\\
&6&5&\dca{6}&\\
\hdashline
&0&0&0&\dca{4}\\
&&&&&\ddots
}.
\]
$[U^{(3)}] = [B_{(1)}]*[L^{(1)}] + $ lower terms, where
\[
B_{(1)} =   \bM{[cc:cc:cc]
\ddots&&&&&\\
\ddots&*&&&\\
\hdashline
&6&*&&\\
&&0&*&\\
\hdashline
&&&6&*&\\
&&&&\ddots&\ddots
},
\quad
L^{(1)} =   \bM{[cc:cc:cc]
\ddots&&&&&\\
&6&&&\\
\hdashline
&6&4&&\\
&\uca{}&5&6&\\
\hdashline
&&&6&4\\
&&&\uca{}&&\ddots
}.
\]
Hence $L^{(1)} = B_{(2)}$ and
\[
[A] =  [B^{(1)}] *[B^{(2)}]* [B^{(3)}] *[B_{(1)}] *[B_{(2)}] + \textup{ lower terms}.
\]
\endexa
By Theorem~\ref{thm:min gen}, Corollary \ref{cor:mA} and \cite[Lemma~24.2.1]{Lu93}, it is now standard to deduce that there exists a canonical basis $\{\{A\}~|~A\in\Tt_d\}$ of $\Sq$ satisfying
\[
\{A\} = [A] + \sum_{A'<_a A} P_{A,A'}[A'],
\]
where $P_{A,A'}\in v^{-1}\NN[v^{-1}]$.
Such a canonical basis for $\Sq$ was known earlier \cite{Lu99} by different constructions (cf. also \cite{DF14}).

\subsection{Modified quantum affine algebra}
Let $\~{\Tt}$ be the set of $\ZZ\times\ZZ$  matrices over $\ZZ$ in which each element $A = (a_{ij})$ satisfies
 $a_{ij} = a_{n+i,n+j}$ for all $i,j$ and $a_{ij} \in \NN$ for all $ i\neq j$.
Let $\dKq$ be the free $\ZZ[v,v^{-1}]$-module with basis $\{ [A] ~|~A\in \~{\Tt}\}$.
The so-called stabilization property for affine $\mathfrak{gl}_n$ case, which is derived from the multiplication formulas (cf. \cite{DF15}) and any monomial basis (e.g. the one in \eqref{eq:mono} and/or any one in \cite[Proposition~6.2]{DF14}), was given in \cite[Proposition 6.3]{DF14}. With this property, the BLM stabilization procedure gives $\dKq$ a unique (non-unital) associative $\ZZ[v,v^{-1}]$-algebra structure with a bar-involution, which is isomorphic to the modified quantum affine $\fgl_n$ (cf. \cite{DF14}).

There is an algorithm (almost identical to Algorithm \ref{alg:mono}) which provides a monomial basis for $\dUq(\^{\fgl}_n)$.
The construction of monomial basis in turn implies a canonical basis for $\dUq(\^{\fgl}_n)$, as already shown in \cite{DF14}.
\rmk
In a recent work \cite{LW15} Li and Wang showed that the canonical basis of $\dUq(\fgl_n)$ does not have positive structure constants in general,
and consequently neither does the canonical basis for $\dUq(\^{\fgl}_n)$.
Hence the \textit{stably-canonical basis} could be a more suitable name for the basis of $\dUq(\^{\fgl}_n)$ arising from the BLM construction. Moreover, Fu and Shoji showed in \cite{FS14} that the positivity property still holds for $\dot{\textbf{U}}_q(\hat{\mathfrak{sl}}_n)$.
\endrmk

\begin{thebibliography}{DDPW08}\frenchspacing

\bibitem[BB05]{GTM231}
A.~Bj{\"o}rner and F.~Brenti, \emph{Combinatorics of {C}oxeter groups},
  Graduate Texts in Mathematics, vol. 231, Springer, New York, 2005.

\bibitem[BKLW14]{BKLW14}
H.~Bao, J.~Kujawa, Y.~Li, and W.~Wang, \emph{Geometric {S}chur duality of
  classical type}, (with appendix by Bao, Li and Wang), Transformation Groups, to appear, arXiv:1404.4000v3 (2014).

\bibitem[BLM90]{BLM90}
A.~Beilinson, G.~Lusztig, and R.~MacPherson, \emph{A geometric setting for the
  quantum deformation of {${\rm GL}_n$}}, Duke Math. J. \textbf{61} (1990),
  no.~2, 655--677.

\bibitem[Cur85]{Cur85}
Charles~W. Curtis, \emph{On {L}usztig's isomorphism theorem for {H}ecke
  algebras}, J. Algebra \textbf{92} (1985), no.~2, 348--365.

\bibitem[DDF12]{DDF12}
B.~Deng, J.~Du, and Q.~Fu, \emph{A double {H}all algebra approach to affine quantum {S}chur-Weyl theory}, London Mathematival Society Lecture Note Series, vol. 401, Cambridge University Press, 2012.

\bibitem[DDPW08]{DDPW08}
B.~Deng, J.~Du, B.~Parshall, and J.~Wang, \emph{Finite dimensional algebras and
  quantum groups}, Mathematical Surveys and Monographs, vol. 150, American
  Mathematical Society, Providence, RI, 2008.

\bibitem[DDX07]{DDX07}
B.~Deng, J.~Du, and J.~Xiao, \emph{Generic extensions and canonical bases for
  cyclic quivers}, Canad. J. Math. \textbf{59} (2007), no.~6, 1260--1283.

\bibitem[DF14]{DF14}
J.~Du and Q.~Fu, \emph{The integral quantum loop algebra of $\mathfrak{gl}_n$},
  arXiv:1404.5679 (2014).

\bibitem[DF15]{DF15}
J.~Du and Q.~Fu, \emph{Quantum affine $\mathfrak{gl}_n$ via {H}ecke algebras},
Adv. Math., {\bf 282} (2015), 23--46.

\bibitem[Du92]{Du92}
J.~Du, \emph{Kazhdan-{L}usztig bases and isomorphism theorems for {$q$}-{S}chur
  algebras}, Kazhdan-{L}usztig theory and related topics ({C}hicago, {IL},
  1989), Contemp. Math., vol. 139, Amer. Math. Soc., Providence, RI, 1992,
  pp.~121--140.

\bibitem[FL14]{FL14}
Z.~Fan and Y.~Li, \emph{Geometric {S}chur duality of classical type, II},
  Trans. Amer. Math. Soc., Series B. {\bf 2} (2015), 51--92.

\bibitem[FL$^3$W16]{FLLLW16}
Z.~Fan, C.~Lai, Y.~Li, L.~Luo, and W.~Wang,
\emph{Affine Hecke algebras and quantum symmetric pairs},
  arXiv:1609.06199 (2016).

\bibitem[FS14]{FS14}
Q.~Fu and T.~Shoji, \emph{Positivity properties for canonical bases of modified quantum affine $\mathfrak{sl}_n$}, arXiv:1407.4228 (2014).

\bibitem[Gre99]{Gr99}
R.~M. Green, \emph{The affine {$q$}-{S}chur algebra}, J. Algebra \textbf{215}
  (1999), no.~2, 379--411.

\bibitem[GV93]{GV93}
V.~Ginzburg and E.~Vasserot, \emph{Langlands reciprocity for affine quantum
  groups of type {$A_n$}}, Internat. Math. Res. Notices (1993), no.~3, 67--85.

\bibitem[KL79]{KL79}
D.~Kazhdan and G.~Lusztig, \emph{Representations of {C}oxeter groups and
  {H}ecke algebras}, Invent. Math. \textbf{53} (1979), no.~2, 165--184.

\bibitem[Lus93]{Lu93}
G.~Lusztig, \emph{Introduction to quantum groups}, Progress in Mathematics,
  vol. 110, Birkh\"auser Boston, Inc., Boston, MA, 1993.

\bibitem[Lus99]{Lu99}
G.~Lusztig, \emph{Aperiodicity in quantum affine {$\mathfrak{gl}_n$}}, Asian J.
  Math. \textbf{3} (1999), no.~1, 147--177.

\bibitem[LW15]{LW15}
Y.~Li and W.~Wang, \emph{Positivity vs negativity of canonical bases},
  arXiv:1501.00688v2 (2015).

\bibitem[VV99]{VV99}
M.~Varagnolo and E.~Vasserot, \emph{On the decomposition matrices of the
  quantized {S}chur algebra}, Duke Math. J. \textbf{100} (1999), no.~2,
  267--297.

\end{thebibliography}

\end{document}